\documentclass[JEDP]{cedram-sem}
\usepackage{amsfonts}
\usepackage{amsthm}
\theoremstyle{plain}
\newtheorem{theorem}{Theorem}
\newtheorem{proposition}[theorem]{Proposition}
\newtheorem{lemma}[theorem]{Lemma}
\newtheorem{corollary}[theorem]{Corollary}
\theoremstyle{definition}

\newtheorem{remark}[theorem]{Remark}%

\newcommand{\R}{{\mathbb R}}

\newcommand{\OO}{{\mathcal O}}
\newcommand{\dd}{\,{\rm d}}
\newcommand{\D}{{\rm d}}
\renewcommand{\div}{\mathop{\mathrm{div}}}
\newcommand{\curl}{\mathop{\mathrm{curl}}}

\newcommand{\meas}{\mathop{\mathrm{meas}}}
\newcommand{\loc}{\mathrm{loc}}
\newcommand{\one}{\mathbf{1}}
\renewcommand{\:}{\thinspace :}

\title 
[The Navier-Stokes Equation in an Exterior Domain]
{Long-Time Asymptotics for the Navier-Stokes Equation
in a Two-Dimensional Exterior Domain}

\alttitle{Comportement asymptotique en temps des solutions
de l'\'equation \\de Navier-Stokes dans un domaine ext\'erieur 
du plan}

\author
[Th. \lastname{Gallay}]
{\firstname{Thierry} \lastname{Gallay}}

\address{
Institut Fourier\\
Universit\'e de Grenoble 1\\
100, rue des Maths, B.P. 74\\
38402 Saint-Martin-d'H\`eres, France}

\thanks{The author is partially supported by the ANR project PREFERED}
\email{Thierry.Gallay@ujf-grenoble.fr}

  
\altkeywords{Equation de Navier-Stokes, comportement asymptotique, 
domaine ext\'erieur}


\begin{document}

\begin{abstract}
We study the long-time behavior of infinite-energy solutions to the
incompressible Navier-Stokes equations in a two-dimensional exterior
domain, with no-slip boundary conditions. The initial data we
consider are finite-energy perturbations of a smooth vortex with
small circulation at infinity, but are otherwise arbitrarily
large. Using a logarithmic energy estimate and some interpolation
arguments, we prove that the solution approaches a self-similar 
Oseen vortex as $t \to \infty$. This result was obtained in 
collaboration with Y.~Maekawa (Kobe University). 
\end{abstract}

\begin{altabstract}
Nous \'etudions le comportement asymptotique en temps des solutions
de l'\'equa\-tion de Navier-Stokes incompressible dans un domaine 
ext\'erieur du plan, avec condition de non-glissement \`a la fronti\`ere.
Les donn\'ees initiales que nous consid\'erons sont des perturbations
d'\'energie finie d'un tourbillon r\'egulier dont la circulation \`a 
l'infini est petite, mais nous n'imposons aucune autre restriction 
\`a leur taille. En utilisant une estimation d'\'energie logarithmique
et des arguments d'interpolation, nous montrons que la solution
converge lorsque $t \to \infty$ vers un tourbillon d'Oseen 
auto\-similaire. Ce r\'esultat a \'et\'e obtenu en collaboration avec
Y.~Maekawa (Universit\'e de Kobe). 
\end{altabstract}

\maketitle

\section{Introduction}\label{sec1}

We consider the free motion of an incompressible viscous fluid
in a two-dimensional exterior domain $\Omega = \R^2 \setminus K$, 
where $K \subset \R^2$ is a compact obstacle with a smooth boundary. 
We do not assume that $K$ is connected, hence we include the case
where the fluid moves around a finite collection of obstacles. 
As for the boundary conditions, we suppose that the velocity 
of the fluid vanishes on $\partial\Omega$ and decays to zero 
at infinity. The evolution of our system is thus governed by the 
Navier-Stokes equations
\begin{equation}\label{NS}
  \left\{\begin{array}{llll}
  \partial_t u + (u\cdot\nabla)u \,=\, \Delta u - \nabla p~, 
  \quad \div u \,=\, 0~, & \quad & \hbox{for~} x \in \Omega\,, \quad
  &t > 0\,, \\
  u(x,t) \,=\, 0~, & & \hbox{for~} x \in \partial\Omega\,,
  &t > 0\,, \\
  u(x,0) \,=\, u_0(x)~, & & \hbox{for~} x \in \Omega\,,&
  \end{array}\right.
\end{equation}
where $u(x,t) \in \R^2$ and $p(x,t) \in \R$ denote, respectively, the
velocity and the pressure of the fluid at a space-time point $(x,t)
\in \Omega \times \R_+$. As can be seen from the first equation in
\eqref{NS}, we assume that the fluid density is constant and that the
kinematic viscosity is equal to $1$. Since \eqref{NS} includes no
forcing, the motion of the fluid originates entirely from the initial 
data $u_0: \Omega \to \R^2$, which we assume to be divergence-free 
and tangent to the boundary on $\partial\Omega$.

The behavior of the solutions of \eqref{NS} depends in a crucial way on 
the decay rate of the velocity field $u(x,t)$ as $|x| \to \infty$.
If the initial velocity $u_0$ belongs to the energy space
\[
  L^2_\sigma(\Omega) \,=\, \Bigl\{u \in L^2(\Omega)^2 \,\Big|\, 
  \div u = 0 \hbox{ in }\Omega\,,~ u\cdot n = 0 \hbox{ on }
  \partial\Omega\Bigr\}~, 
\]
where $n$ denotes the interior unit normal on $\partial \Omega$, we 
have the following classical result\:

\begin{theorem}\label{thm1}
For all initial data $u_0 \in L^2_\sigma(\Omega)$, Eq.~\eqref{NS} has 
a unique global solution
\[
  u \in C^0([0,\infty), L^2_\sigma(\Omega)) \cap C^1((0,\infty), 
  L^2_\sigma(\Omega)) \cap C^0((0,\infty), H^1_0(\Omega)^2\cap 
  H^2(\Omega)^2)~,
\]
which satisfies for all $t \ge 0$ the energy equality\:
\begin{equation}\label{Eeq}
  \frac12 \|u(\cdot,t)\|_{L^2(\Omega)}^2 + \int_0^t \|\nabla u(\cdot,
  s)\|_{L^2(\Omega)}^2\dd s \,=\, \frac12 \|u_0\|_{L^2(\Omega)}^2~.
\end{equation}
\end{theorem}

Global well-posedness for the Navier-Stokes equations was first
established by Leray \cite{Le1} in the particular case where
$\Omega = \R^2$. When $\Omega \subset \R^2$ is bounded, the
first results also go back to Leray \cite{Le2}, but global existence
of large solutions was shown only later by Ladyzhenskaia \cite{La},
see also \cite{LP,KF,FK}. To prove Theorem~\ref{thm1}, one can
use a regularization or a discretization procedure to construct global 
{\em weak} solutions of \eqref{NS} which satisfy the energy inequality, 
and then prove that these solutions are unique and have the desired 
regularity. Alternatively, one can construct local {\em mild} solutions 
by transforming \eqref{NS} into an integral equation and solving it by 
a fixed point argument, and then use the energy equality \eqref{Eeq} 
to show that all solutions can be extended to the whole time interval
$[0,\infty)$. Although most of the literature is devoted to
the situation where $\Omega$ is either a bounded domain or the whole
plane $\R^2$, the case of an exterior domain can be treated without
essential modifications, see e.g. \cite{KOM}.

It follows from \eqref{Eeq} that the kinetic energy $E(t) = \frac12
\|u(\cdot,t)\|_{L^2(\Omega)}^2$ is nonincreasing in time, and a result
of Masuda \cite{Ma} shows that $E(t)$ converges to zero as $t \to
\infty$. Moreover, under additional assumptions on the initial data,
it is possible to specify a decay rate in time. For instance, if 
$u_0 \in L^2_\sigma (\Omega) \cap L^q(\Omega)^2$ for some $q \in (1,2)$, 
the solution of \eqref{NS} lies in the same space for all $t > 0$ 
and
\begin{equation}\label{Edecay}
  \|u(\cdot,t)\|_{L^2(\Omega)} \,=\, o\Bigl(t^{\frac12 - \frac1q}\Bigr)
  \qquad \hbox{as} \quad t \to \infty~,
\end{equation}
see \cite{BM,KOA,BJ}. It is interesting to notice that \eqref{Edecay}
fails in the limiting case $q = 1$. Indeed, if $u_0 \in L^2_\sigma
(\Omega) \cap L^1(\Omega)^2$, then in general the velocity field
$u(x,t)$ decays like $|x|^{-2}$ as $|x| \to \infty$, so that 
$u(\cdot,t) \notin L^1(\Omega)^2$ for $t > 0$.  

As an aside, we mention that this loss of spatial decay is related to
the net force $F$ exerted by the fluid on the obstacle $K$. To see
this, we first observe that any velocity field $u \in L^2_\sigma
(\Omega) \cap L^1(\Omega)^2$ satisfies $\int_\Omega u \dd x = 0$. 
Indeed, if $u$ is smooth and compactly supported, then using Gauss' 
theorem and the fact that $u \cdot n = 0$ on $\partial\Omega$ we find
\[
  \int_\Omega u_j \dd x \,=\, \int_\Omega (u\cdot\nabla)x_j \dd x 
  \,=\, \int_\Omega \div(u\,x_j) \dd x \,=\, 0~, \qquad 
  \hbox{for }j = 1,2~.
\]
The general easily case follows by a density argument \cite{Koz}. 
Now, if $u \in C^1([0,T],L^1(\Omega)^2)$ is a solution of the 
Navier-Stokes equation \eqref{NS}, then
\[
  0 \,=\, \frac{\D}{\D t} \int_\Omega u\dd x \,=\, \int_\Omega 
  \Bigl(\Delta u - \nabla p - (u\cdot\nabla)u\Bigr)\dd x 
  \,=\, -\int_{\partial\Omega} (Tn)\dd\sigma \,=\, -F~,
\]
because $\Delta u_i -\partial_i p = \partial_j T_{ij}$ where $T_{ij} = 
\partial_i u_j + \partial_j u_i - p\delta_{ij}$ is the stress tensor
(we recall that all physical parameters have been normalized to $1$).
The formal calculation above can be made rigorous \cite{Koz} and shows
that, no matter how localized the initial data may be, the velocity
field $u(\cdot,t)$ does not stay integrable for positive times, unless
the net force $F$ vanishes identically. Of course this is not the case
in general, but in highly symmetric situations it is possible to
construct solutions of \eqref{NS} for which $F \equiv 0$, and which
decay faster as $t \to \infty$ than what is indicated in
\eqref{Edecay}, see \cite{HM1,HM2}.

Much less is known about the solutions of \eqref{NS} if the initial
data $u_0$ are not square integrable. Although the physical relevance
of infinite-energy solutions can be questioned, we believe that such
solutions naturally occur when studying the dynamics of \eqref{NS}
in a two-dimensional exterior domain $\Omega$. One way to realize that
is to consider the relation between the velocity field $u$ and the
associated vorticity $\omega = \partial_1 u_2 -\partial_2 u_1$. Given
$p \in [1,2)$, we denote
\[
  \dot W_{0,\sigma}^{1,p}(\Omega) \,=\, \Bigl\{u \in L^\frac{2 p}{2-p}
  (\Omega )^2 ~\Big |~ \nabla u \in L^p(\Omega)^4\,,~\div u = 0 
  ~\hbox{in}~\Omega\,,~ u=0~\hbox{on}~\partial\Omega\Bigr\}~.
\]
In other words $\dot W^{1,p}_{0,\sigma}(\Omega)$ is the completion with 
respect to the norm $u \mapsto \|\nabla u\|_{L^p}$ of the space of all 
smooth, divergence-free vector fields with compact support in $\Omega$, 
see \cite{Gal}. We then have the following result\:

\begin{lemma}\label{BS0}
If $u \in \dot W^{1,p}_{0,\sigma}(\Omega)$ for some $p \in [1,2)$, and 
if $\omega = \partial_1 u_2 - \partial_2 u_1$, then
\begin{equation}\label{BS}
  u(x) \,=\, \frac{1}{2\pi}\int_{\Omega} \frac{(x-y)^\perp}{
  |x-y|^2}\,\omega(y)\dd y~, 
\end{equation}
for almost every $x \in \Omega$. Here $x^\perp = (-x_2,x_1)$ and 
$|x|^2 = x_1^2 + x_2^2$ if $x = (x_1,x_2) \in \R^2$. 
\end{lemma}

The proof of Lemma~\ref{BS0} is very simple\: if $\bar u : \R^2 \to 
\R^2$ denotes the extension of $u$ by zero outside $\Omega$, then
$\bar u \in L^{2 p/(2-p)}(\R^2)^2$, $\nabla \bar u \in L^p(\R^2)^4$, 
and $\div \bar u = 0$. Moreover $\bar \omega = \partial_1 \bar u_2 - 
\partial_2 \bar u_1$ is the extension of $\omega$ by zero outside
$\Omega$. Thus $\bar u$ can be expressed in terms of $\bar \omega$
using the classical Biot-Savart law in $\R^2$, and restricting that
relation to $\Omega$ we obtain \eqref{BS}. We emphasize that the
representation \eqref{BS} is only valid if $\omega$ is the curl of a
divergence-free velocity field $u$ which {\em vanishes on}
$\partial\Omega$. In contrast, if $\omega : \Omega \to \R$ is an
arbitrary smooth function with compact support, the velocity field
{\em defined by} \eqref{BS} does not even satisfy $u \cdot n = 0$ on
$\partial \Omega\,$!

We now assume that the vorticity distribution $\omega$ is sufficiently
localized so that $\omega \in L^1(\Omega)$, and we define the 
{\em total circulation}
\[
  \alpha \,=\, \int_\Omega \omega(x) \dd x \,=\, \lim_{R \to \infty}
  \oint_{|x| = R} u_1 \dd x_1 + u_2 \dd x_2~,
\]
where the second equality follows from Green's theorem, since $\omega
= \partial_1 u_2 - \partial_2 u_1$ and $u$ vanishes on $\partial
\Omega$. Using the vorticity formulation of the Navier-Stokes
equation \eqref{NS}, it is not difficult to verify that the total
circulation is a conserved quantity. But it follows from \eqref{BS}
that
\begin{equation}\label{uinfty}
  u(x) ~\sim~ \frac{\alpha}{2\pi}\,\frac{x^\perp}{|x|^2}~, \qquad 
  \hbox{as}\quad |x| \to \infty~,
\end{equation}
hence $u \notin L^2(\Omega)^2$ as soon as $\alpha \neq 0$. This 
shows that finite-energy solutions of \eqref{NS} necessarily have 
zero total circulation. On the other hand, in many important 
examples of two-dimensional flows such as vortex patches, vortex 
sheets, or point vortices, the vorticity distribution typically
has a constant sign, hence the total circulation is necessarily 
nonzero. In our opinion, it is thus important to enlarge the class 
of admissible solutions of \eqref{NS}, so as to allow for velocity 
fields which decay like $|x|^{-1}$ as $|x| \to \infty$. 

A possible framework for the study of infinite-energy solutions 
of the Navier-Stokes equation \eqref{NS} is the {\em weak energy
space} 
\[
  L^{2,\infty}_\sigma(\Omega) \,=\,  \Bigl\{u \in L^{2,\infty}
  (\Omega)^2\,\Big|\, \div u \,=\, 0 \hbox{~in~} \Omega\,,~ 
  u\cdot n = 0 \hbox{~on~}\partial\Omega\Bigr\}~,
\]
where $L^{2,\infty}(\Omega)$ is the weak $L^2$ space on $\Omega$, 
see \cite{BL}. We recall that
\begin{equation}\label{weaknorm}
  \|u\|_{L^{2,\infty}(\Omega)} \,\approx\, \sup_{\lambda > 0} \lambda \Bigl(
  \meas\{x \in \Omega \,|\, |u(x)| > \lambda\}\Bigr)^{1/2}~,
\end{equation}
in the sense that the norm $\|u\|_{L^{2,\infty}}$ is equivalent to the
quantity in the right-hand side of \eqref{weaknorm}.  Clearly
$L^2_\sigma(\Omega) \hookrightarrow L^{2,\infty}_\sigma(\Omega)$, but
the weak energy space is large enough to include velocity fields which
decay slowly at infinity, as in \eqref{uinfty}.  Concerning the
solvability of \eqref{NS} in $L^{2,\infty}_\sigma(\Omega)$, the
following general result was obtained by Kozono and Yamazaki\:

\begin{theorem}\label{thm3} \cite{KY} 
There exists $\epsilon > 0$ such that, for all initial data $u_0 
\in L^{2,\infty}_\sigma(\Omega)$ satisfying
\begin{equation}\label{limsing}
  \limsup_{\lambda \to +\infty} \lambda \Bigl(\meas\{x \in 
  \Omega \,|\, |u_0(x)| > \lambda\}\Bigr)^{1/2} \,\le\, \epsilon~,
\end{equation}
Eq.~\eqref{NS} has a unique global solution such that, for all 
$T > 0$, 
\[
  \sup_{0 < t < T} \|u(\cdot,t)\|_{L^{2,\infty}(\Omega)} + 
  \sup_{0 < t < T} t^{1/4}\|u(\cdot,t)\|_{L^4(\Omega)} \,<\, \infty~, 
\]
and such that $u(\cdot,t) \to u_0$ as $t \to 0$ in the weak-$*$ 
topology of $ L^{2,\infty}_\sigma(\Omega)$.
\end{theorem}

Theorem~\ref{thm3} shows that the Cauchy problem for the Navier-Stokes
equations \eqref{NS} is globally well-posed in the weak energy space
$L^{2,\infty}_\sigma(\Omega)$, provided that the local singularity of
the initial data $u_0$ is sufficiently small, in the sense of
\eqref{limsing}. To illustrate the meaning of this smallness
condition, we consider the simple situation where the initial flow is
just a point vortex of circulation $\alpha \in \R$ located at $x_0 \in
\Omega$. The initial vorticity is thus given by $\omega_0(x) = \alpha
\delta(x-x_0)$, and using the classical Biot-Savart law in the
exterior domain $\Omega$ (see e.g. \cite{ILN}) it is easy to verify
that the corresponding velocity field $u_0$ lies in
$L^{2,\infty}_\sigma(\Omega)$, is smooth in $\Omega \setminus
\{x_0\}$, and satisfies
\[
  u_0(x) \,\approx\, \frac{\alpha}{2\pi} \,\frac{(x-x_0)^\perp} 
  {|x-x_0|^2}\qquad \hbox{as}\quad x \to x_0~,
\]
so that \eqref{limsing} is fulfilled if and only if $|\alpha| \le
\sqrt{4\pi}\epsilon$. This example shows that, if the initial
vorticity $\omega_0$ is a finite measure, condition \eqref{limsing}
implies a restriction on the size of the atomic part of $\omega_0$.
Such a restriction also arises in the analysis of the two-dimensional
vorticity equation in the whole space $\R^2$, see \cite{GMO}, but 
in that particular case the uniqueness of the solution can be 
established when the initial vorticity is an arbitrary finite measure 
\cite{GG,BedMas}. 

Although Theorem~\ref{thm3} provides the existence of a large class of
infinite-energy solutions, very little is known about the asymptotic
behavior of these solutions as $t \to \infty$. In fact, we do not even
know whether they stay bounded in the weak energy space
$L^{2,\infty}_\sigma(\Omega)$, because we are lacking a priori
estimates. Indeed, if $u_0 \notin L^2_\sigma(\Omega)$ the energy
equality \eqref{Eeq} does not make sense, and because of the no-slip
boundary condition on $\partial\Omega$ it is quite difficult to obtain
estimates on the vorticity distribution if $\Omega \neq \R^2$. In the
rest of this paper, however, we consider a particular class of
infinite-energy solutions of the Navier-Stokes equations \eqref{NS},
for which the asymptotic behavior in time can be accurately described.

\section{Main Results}\label{sec2}

In the particular case where $\Omega = \R^2$, the Navier-Stokes 
equations \eqref{NS} have a family of self-similar solutions of the 
form $u(x,t) \,=\, \alpha \Theta(x,t)$, $p(x,t) \,=\, \alpha^2 
\Pi(x,t)$, where $\alpha \in \R$ is a free parameter (the total 
circulation) and
\begin{equation}\label{Thetadef}
  \Theta(x,t) \,=\, \frac{1}{2\pi}\,\frac{x^\perp}{|x|^2}
   \Bigl(1 - e^{-\frac{|x|^2}{4(1+t)}}\Bigr)~, \qquad
  \nabla \Pi(x,t) \,=\, \frac{x}{|x|^2}|\Theta(x,t)|^2~.
\end{equation}
These solutions are usually called the {\it Lamb-Oseen vortices}. 
If $u(x,t) \,=\, \alpha \Theta(x,t)$, the corresponding vorticity 
distribution is $\omega(x,t) = \alpha \Xi(x,t)$, where 
\begin{equation}\label{Xidef}
  \Xi(x,t) \,=\, \partial_1 \Theta_2(x,t) - \partial_2 \Theta_1(x,t) 
  \,=\, \frac{1}{4\pi(1+t)}\,e^{-\frac{|x|^2}{4(1+t)}}~.
\end{equation}
Note that $\Xi(x,t) > 0$ and $\int_{\R^2} \Xi(x,t)\dd x = 1$ for
all $t \ge 0$. Oseen vortices play an important role in the dynamics
of the Navier-Stokes equations in $\R^2$. In particular, we have 
the following result\:

\begin{theorem}\label{thm4} \cite{GW2}. 
For all initial data $u_0 \in L^{2,\infty}_\sigma(\R^2)$ such that 
the vorticity distribution $\omega_0$ is integrable, the solution
of the Navier-Stokes equation in $\R^2$ satisfies
\[
  \int_{\R^2} |\omega(x,t) - \alpha \Xi(x,t)|\dd x 
  \,\xrightarrow[t\to\infty]{}\, 0~, \qquad \hbox{where}
  \quad \alpha \,=\, \int_{\R^2} \omega_0 \dd x~.
\]
\end{theorem}

In other words, Oseen vortices describe the leading order asymptotics
of all solutions of the Navier-Stokes equations in $\R^2$ with 
integrable vorticity distribution and nonzero total circulation. 

In the case of an exterior domain $\Omega = \R^2\setminus K$,
approximate Oseen vortices can be constructed in the following way. Let
$\chi : \R^2 \to [0,1]$ be a smooth, radially symmetric cut-off
function such that $\chi$ is nondecreasing along rays, $\chi = 0$ on a
neighborhood of $K$,  and $\chi(x) = 1$ when $|x|$ is sufficiently large.
The {\em truncated Oseen vortex} with unit circulation is defined as 
follows\:
\begin{equation}\label{uchidef}
  u^\chi(x,t) \,=\,  \chi(x)\Theta(x,t) \,=\, \frac{1}{2\pi}
  \,\frac{x^\perp}{|x|^2} \Bigl(1 - e^{-\frac{|x|^2}{4(1+t)}}\Bigr)\chi(x)~.
\end{equation}
For any $t \ge 0$, it is clear that $u^\chi(\cdot,t)$ is a smooth 
divergence-free vector field which vanishes in a neighborhood of 
$K$. The corresponding vorticity distribution $\omega^\chi = 
\partial_1 u^\chi_2 - \partial_2 u^\chi_1$ has the explicit
expression
\begin{equation}\label{omegachi}
  \omega^\chi(x,t) \,=\, \chi(x)\Xi(x,t) + \frac{1}{2\pi}\,
  \frac{1}{|x|^2}\Bigl(1 - e^{-\frac{|x|^2}{4(1+t)}}\Bigr)
  x\cdot\nabla\chi(x)~,
\end{equation}
where $\Xi(x,t)$ is defined in \eqref{Xidef}. In particular 
$\omega^\chi(x,t) \ge 0$ and $\int_{\R^2}\omega^\chi(x,t) \dd x = 1$
for all $t \ge 0$. Moreover, a direct calculation shows that
\begin{equation}\label{uchirel}
  (u^\chi \cdot \nabla)u^\chi \,=\, \frac12 \nabla|u^\chi|^2 + 
  (u^\chi)^\perp \omega^\chi \,=\, -\frac{x}{|x|^2}|u^\chi|^2~,
\end{equation}
hence there exists a radially symmetric function $p^\chi(x,t)$ such
that $-\nabla p^\chi = (u^\chi \cdot \nabla)u^\chi$. 

Now, given $\alpha \in \R$, we consider solutions of \eqref{NS} of 
the particular form
\begin{equation}\label{udec}
  u(x,t) \,=\, \alpha u^\chi(x,t) + v(x,t)~, \qquad 
  p(x,t) \,=\, \alpha^2 p^\chi(x,t) + q(x,t)~,
\end{equation}
where $u^\chi(x,t)$ is the truncated Oseen vortex defined in 
\eqref{uchidef}, and $v(x,t)$ is a {\em finite-energy} perturbation. 
In this situation, we expect that $v(\cdot,t)$ converges to zero 
in energy norm as $t \to \infty$, so that the long-time behavior 
of $u(\cdot,t)$ is described, to leading order, by the Oseen 
vortex $\alpha \Theta(\cdot,t)$. Our main result, which was obtained
in collaboration with Y.~Maekawa, shows that this is indeed the case, 
provided the total circulation $\alpha$ is sufficiently small. 

\begin{theorem}\label{main} \cite{GM}
Fix $q \in (1,2)$, and let $\mu = 1/q - 1/2$. There exists a constant 
$\epsilon = \epsilon(q) > 0$ such that, for any smooth exterior 
domain $\Omega \subset \R^2$ and for all initial data of the form 
$u_0 = \alpha u^\chi(\cdot,0) + v_0$ with $|\alpha| \le \epsilon$ 
and $v_0 \in L^2_\sigma(\Omega) \cap L^q(\Omega)^2$, the solution of the 
Navier-Stokes equations \eqref{NS} satisfies
\begin{equation}\label{conv}
  \|u(\cdot,t) -\alpha\Theta(\cdot,t)\|_{L^2(\Omega)} + 
  t^{1/2} \|\nabla u(\cdot,t) -\alpha\nabla\Theta(\cdot,t)
  \|_{L^2(\Omega)} \,=\, \OO(t^{-\mu})~,
\end{equation}
as $t \to \infty$.  
\end{theorem}

To understand the scope and the limitations of this statement, a few
comments are in order. First of all, Theorem~\ref{main} is a {\em
global stability result} for Oseen vortices with small circulation
at infinity, because we do not impose any restriction on the size of
the perturbation $v_0 \in L^2_\sigma(\Omega) \cap L^q(\Omega)^2$. In 
the particular case where $\alpha = 0$, there is no vortex at all and
we just recover the asymptotics \eqref{Edecay} with $\OO(t^{-\mu})$ 
instead of $o(t^{-\mu})$ in the right-hand side. Also, in the 
simple situation where $\Omega = \R^2$, our result is comparable 
to that of Carpio \cite{Car}, although the proof is very different. 

The main limitation of Theorem~\ref{main} is of course the restriction
on the size of the circulation $\alpha$, which we believe is purely
technical. In this respect, the fact that $\epsilon(q)$ can be taken
independent of the domain $\Omega$ is quite significant, because we
know that there is no restriction on the circulation in the particular
case where $\Omega = \R^2$, see Theorem~\ref{thm4}. Obviously
$\epsilon(q)$ is a decreasing function of $q$, and the proof shows
that $\epsilon(q) = \OO(\sqrt{2-q})$ as $q \to 2$. Thus the limiting
case $q = 2$ is not included, which means that we are not able to
control arbitrary finite-energy perturbations of the Oseen vortex
(see however \cite{IKL} for a partial result in that direction).  On
the other hand the limit of $\epsilon(q)$ as $q \to 1$ can be
estimated and is found to be approximately $\epsilon_* = 5.306$, see
\cite{GM}.

We also mention that the decomposition $u_0 = \alpha u^\chi(\cdot,0) 
+ v_0$ of the initial data is automatically satisfied if we assume 
that the initial vorticity is sufficiently localized. More precisely, 
we have the following auxiliary result, which follows quite 
easily from Lemma~\ref{BS0}. 

\begin{proposition}\label{vort1} \cite{GM}
Given $q\in (1, 2)$, assume that $u_0 \in \dot W_{0,\sigma}^{1,p}(\Omega)$ 
for some $p\in [1,2)$ and that the associated vorticity $\omega_0 = 
\curl u_0$ satisfies 
\begin{equation}\label{omloc}
  \int_{\Omega} (1 + |x|^2)^m |\omega_0 (x)|^2 \dd x \,<\,\infty~, 
\end{equation}
for some $m > 2/q$. If we denote $\alpha = \int_{\Omega} \omega_0 (x) 
\dd x$, then $u_0 = \alpha u^\chi(\cdot,0) + v_0$ for some $v_0 \in 
L^2_\sigma(\Omega) \cap L^q(\Omega)^2$. In particular, if $|\alpha| \le 
\epsilon$, the conclusion of Theorem \ref{main} holds.
\end{proposition}
 
In view of Proposition~\ref{vort1}, it would be nice to extend the
conclusion Theorem~\ref{main} so as to include a convergence result
for the vorticity distribution in the critical space
$L^1(\Omega)$. This is not immediately obvious, because the classical
$L^p$--$L^q$ estimates for the Stokes semigroup in an exterior domain
do not include the limiting case $p = 1$, see \cite{DS1,DS2}.
However, combining Theorem~\ref{main} with a relatively standard
estimate, which shows that the $L^1$ norm of the vorticity cannot leak 
to infinity, we obtain the following result which is the main original
contribution of the present paper\:

\begin{proposition}\label{vort2}
Under the assumptions of Theorem~\ref{main}, if we suppose in addition 
that \eqref{omloc} holds for some $m \ge 2/q$, then the vorticity 
$\omega = \curl u$ satisfies
\begin{equation}\label{convom}
  \int_\Omega |\omega(x,t) - \alpha \Xi(x,t)|\dd x \,=\, 
  \OO(t^{-\mu}\log t)~, \qquad \hbox{as}\quad t \to \infty~,
\end{equation}
where $\Xi(x,t)$ is defined in \eqref{Xidef}. 
\end{proposition}

In the rest of this paper, we give a simplified proof of
Theorem~\ref{main}, which does not yield the optimal conclusion.  In
particular, we shall find a suboptimal convergence rate in
\eqref{conv}, and our limitation on the size of the circulation will
(a priori) depend on the domain $\Omega$. We refer the reader to
\cite{GM} for a complete proof, including all details. In the last
section, we briefly show how Proposition~\ref{vort2} follows from
Theorem~\ref{main}, using some additional information on the vorticity
near infinity.

\section{Energy estimates}\label{sec3}

Given $\alpha \in \R$ we consider solutions of \eqref{NS} of the form
\eqref{udec}. The perturbation $v(x,t)$ vanishes on the boundary 
$\partial\Omega$ and satisfies the equation
\begin{equation}\label{veq}
  \partial_t v + \alpha (u^\chi\cdot\nabla)v + \alpha (v\cdot\nabla)u^\chi
  + (v\cdot\nabla)v \,=\, \Delta v + \alpha R^\chi - \nabla q~,  
  \qquad \div v \,=\, 0~,
\end{equation}
where $R^\chi$ is the remainder term given by \eqref{Rchidef} below.
It is not difficult to verify that the Cauchy problem for equation \eqref{veq}
is globally well-posed in the energy space $L^2_\sigma(\Omega)$. The
goal of this section is to control the long-time evolution of the 
perturbation $v(t) \equiv v(\cdot,t)$ using energy estimates.

First of all, we multiply both sides of \eqref{veq} by $v$ and 
integrate by parts over $\Omega$. Taking into account the no-slip 
boundary condition, we find
\begin{equation}\label{venergy}
  \frac12 \frac{\D}{\D t} \|v(t)\|_{L^2}^2 + \|\nabla v(t)\|_{L^2}^2 
  \,=\, \alpha \langle v(t), R^\chi(t)\rangle - \alpha \langle v(t), 
  (v(t)\cdot\nabla) u^\chi(t)\rangle~,
\end{equation}
where $\langle \cdot,\cdot \rangle$ denotes the usual scalar product 
in $L^2_\sigma(\Omega)$, so that $\|\cdot\|_{L^2} = \langle \cdot\,,\cdot
\rangle^{1/2}$. To estimate the right-hand side of \eqref{venergy}, 
we use the results of Section~\ref{sec6} below. First, in view of 
\eqref{Rchi2}, we have
\[
  |\alpha \langle v(t),R^\chi(t)\rangle| \,\le\, \frac{\kappa_3\, 
  |\alpha|}{1+t} \|\nabla v(t)\|_{L^2} \,\le\, \frac12 
  \|\nabla v(t)\|_{L^2}^2 + \frac{\kappa_3^2\,\alpha^2}{2(1+t)^2}~.
\]
Moreover, applying \eqref{uchi2} with $p = \infty$, we see that 
\[
  |\langle v(t), (v(t)\cdot\nabla) u^\chi(t)\rangle| \,\le\, 
  \frac{b_\infty}{1+t}\,\|v(t)\|_{L^2}^2~.
\]
We thus obtain the energy inequality
\[
  \frac{\D}{\D t} \|v(t)\|_{L^2}^2 + \|\nabla v(t)\|_{L^2}^2
  \,\le\, \frac{2b_\infty |\alpha|}{1+t}\|v(t)\|_{L^2}^2 + 
  \frac{\kappa_3^2\,\alpha^2}{(1+t)^2}~, \qquad t > 0~.
\]
Using Gronwall's lemma, we deduce that
\begin{equation}\label{vpoly}
  \|v(t)\|_{L^2}^2 + \int_{t_0}^t \|\nabla v(s)\|_{L^2}^2\dd s 
  \,\le\, \Bigl(\frac{1+t}{1+t_0}\Bigr)^{2 b_\infty |\alpha|}
  \Bigl(\|v(t_0)\|_{L^2}^2 + \frac{\kappa_3^2\,\alpha^2}{1+t_0} 
  \frac{t-t_0}{1+t}\Bigr)~,
\end{equation}
for $t \ge t_0 \ge 0$. This simple estimate shows that the energy of 
the perturbation $v(x,t)$ grows at most polynomially in time as
$t \to \infty$. Such a conclusion is rather pessimistic, however, 
because by a relatively simple modification of the previous argument 
it is possible to establish a logarithmic bound, which is clearly 
superior for large times.

\begin{proposition}\label{logbound} There exists a constant $K_1 > 0$
such that, for any circulation $\alpha \in \R$ and any $v_0 \in 
L^2_\sigma(\Omega)$, the solution of \eqref{veq} with initial data 
$v_0$ satisfies
\begin{equation}\label{vlog}
  \|v(t)\|_{L^2}^2 + \int_0^t \|\nabla v(s)\|_{L^2}^2\dd s 
  \,\le\, 2^{c|\alpha|}K_1\Bigl(\|v_0\|_{L^2}^2 + \alpha^2
  \log(1+t)\Bigr)~,
\end{equation}
for all $t \ge 0$, where $c = 2b_\infty > 0$.    
\end{proposition}

\begin{proof}
If $t \le 1$ then \eqref{vlog} follows from \eqref{vpoly} with 
$t_0 = 0$, hence we can assume that $t \ge 1$. Given any $\tau \ge 0$, 
we denote
\begin{equation}\label{tildevdef}
  \tilde v(x,t) \,=\, u(x,t) - \alpha u^\chi(x,t+\tau) \,=\, 
  v(x,t) + \alpha \Bigl(u^\chi(x,t)- u^\chi(x,t+\tau)\Bigr)~, 
\end{equation}
for all $x \in \Omega$ and all $t \ge 0$. Then $\tilde v$ satisfies
\eqref{veq} where $u^\chi(x,t)$ and $R^\chi(x,t)$ are replaced by
$u^\chi(x,t+\tau)$ and $R^\chi(x,t+\tau)$, respectively. Proceeding 
exactly as above, we thus obtain the energy estimate
\begin{equation}\label{tildevpoly}
  \|\tilde v(t)\|_{L^2}^2 + \int_0^t \|\nabla \tilde v(s)\|_{L^2}^2
  \dd s \,\le\, \Bigl(\frac{1+t+\tau}{1+\tau}\Bigr)^{c|\alpha|}
  \Bigl(\|\tilde v(0)\|_{L^2}^2 + C \alpha^2\Bigr)~, 
  \qquad t \ge 0~.
\end{equation}
Now, we fix $t \ge 1$ and choose $\tau = t-1$. From \eqref{uchi3}, 
\eqref{tildevdef}, we have
\[
  \|v(t)\|_{L^2}^2 \,\le\, 2 \|\tilde v(t)\|_{L^2}^2 + 
  2 \alpha^2\|u^\chi(t)  - u^\chi(2t-1)\|_{L^2}^2 \,\le\, 
  2 \|\tilde v(t)\|_{L^2}^2 + 2\kappa_1 \alpha^2 \log 2~.
\]
Similarly, using \eqref{uchi4}, we find
\begin{align*}
 \int_0^t \|\nabla v(s)\|_{L^2}^2 \dd s  \,&\le\, 2 \int_0^t 
  \|\nabla \tilde v(s)\|_{L^2}^2 \dd s  + 2\alpha^2 \int_0^t 
  \|\nabla u^\chi(s)- \nabla u^\chi(s+t-1)\|_{L^2}^2 \dd s \\
  \,&\le\, 2 \int_0^t \|\nabla \tilde v(s)\|_{L^2}^2 \dd s 
  + 2\kappa_2 \alpha^2\,\log\frac{1+t}{2}~.
\end{align*}
Thus it follows from \eqref{tildevpoly} (with $\tau = t-1$) that
\begin{equation}\label{logprelim}
  \|v(t)\|_{L^2}^2 + \int_0^t \|\nabla v(s)\|_{L^2}^2 \dd s  
  \,\le\, 2^{c|\alpha|+1}\Bigl(\|\tilde v(0)\|_{L^2}^2 + C\alpha^2
  \Bigr) + \kappa \alpha^2 \log(1+t)~, 
\end{equation}
for all $t > 0$, where $\kappa = 2\max(\kappa_1,\kappa_2)$. 
Finally, we have by \eqref{uchi3}
\[
   \|\tilde v(0)\|_{L^2}^2 \,\le\, 2 \|v_0\|_{L^2}^2 + 
  2 \alpha^2\|u^\chi(0)  - u^\chi(t-1)\|_{L^2}^2 \,\le\, 
  2 \|v_0\|_{L^2}^2 + 2\kappa_1 \alpha^2 \log t~,
\]
hence \eqref{vlog} easily follows from \eqref{logprelim}.
\end{proof}

\begin{remark}\label{logrem}
The logarithmic energy estimate \eqref{vlog} is the main new ingredient
in the proof of Theorem~\ref{main}. To a certain extent, we use it as 
a substitute for the energy equality \eqref{Eeq}, which does not make
sense for the solutions we consider. As is clear from the proof, 
the logarithmic energy estimate relies on the fact that Oseen's
vortex \eqref{Thetadef} has ``nearly finite energy'', in the sense
that the integral defining $\|\Theta(\cdot,t)\|_{L^2}^2$ diverges only
logarithmically at infinity.
\end{remark}

\section{Fractional interpolation}\label{sec4}

Let $P$ be the Leray-Hopf projection in $\Omega$, and $A = -P\Delta$ 
be the Stokes operator, see e.g. \cite{CF}. We recall that $A$ is
self-adjoint and nonnegative in $L^2_\sigma(\Omega)$, so that the 
fractional power $A^\beta$ can be defined for all $\beta > 0$. The 
following result shows that the range of $A^\mu$ contains the (dense) 
subspace $L^2_\sigma(\Omega) \cap L^q(\Omega)^2$.

\begin{lemma}\label{Amu} {\rm \cite{BM,KOA}} Let $q \in (1,2)$ and
$\mu = \frac1q - \frac12$. For all $v \in L^2_\sigma(\Omega) \cap 
L^q(\Omega)^2$, there exists a unique $w \in D(A^\mu) \subset 
L^2_\sigma(\Omega)$ such that $v = A^\mu w$. Moreover, there exists
a constant $C > 0$ (independent of $v$) such that $\|w\|_{L^2(\Omega)}
\le C\|v\|_{L^q(\Omega)}$. 
\end{lemma}

\begin{remark}\label{A-mu}
If $v,w$ are as in the above statement, we denote $w = A^{-\mu}v$. 
Roughly speaking, the proof of Lemma~\ref{Amu} argues as follows. 
By classical Sobolev embedding, we know that the domain of 
$A^\mu$ is contained in $L^2_\sigma(\Omega) \cap L^{q'}(\Omega)^2$, 
where $\frac{1}{q'} = \frac12-\mu = 1 - \frac1q$, and by duality we 
deduce that the range of $A^{-\mu}$ contains $L^2_\sigma(\Omega) \cap
L^q(\Omega)^2$. 
\end{remark}

We go back to the study of the perturbation equation \eqref{veq}, 
which can be written in the equivalent form
\begin{equation}\label{veq2}
  \partial_t v + A v + \alpha P \Bigl((u^\chi\cdot\nabla)v + 
  (v\cdot\nabla)u^\chi\Bigr) + P(v\cdot\nabla)v \,=\, \alpha 
  R^\chi~. 
\end{equation}
So far we only considered solutions in the energy space $L^2_\sigma
(\Omega)$, but now we assume in addition that $v_0 \in L^q(\Omega)^2$,
for some fixed $q \in (1,2)$, and we denote $\mu = \frac1q - \frac12
\in (0,\frac12)$. Then it is not difficult to verify that the solution
$v(t)$ of \eqref{veq2} lies in $L^2_\sigma(\Omega) \cap L^q(\Omega)^2$
for all $t\ge 0$. In particular, invoking Lemma~\ref{Amu}, we can define
$w(t) = A^{-\mu} v(t)$ for all $t \ge 0$. This quantity solves the
modified equation
\begin{equation}\label{weq}
  \partial_t w + A w + \alpha F_\mu(u^\chi,v) + \alpha F_\mu(v,u^\chi) +
  F_\mu(v,v) \,=\, \alpha A^{-\mu}R^\chi~,
\end{equation}
where $F_\mu(u,v)$ is the bilinear term formally defined by
\begin{equation}\label{Fmudef}
  F_\mu(u,v) \,=\, A^{-\mu} P(u\cdot\nabla)v~.
\end{equation}
We refer to \cite[Section 2]{KOA} for a rigorous definition and 
a list of properties of the bilinear map $F_\mu$. Our goal here
is to establish the following estimate\:

\begin{proposition}\label{wenergy}
There exists $K_3 > 0$ and, for all $\alpha \in \R$, there exist positive 
constants $K_2(\alpha)$ and $k(\alpha)$ such that, if $v$ is any 
solution of \eqref{veq2} with initial data $v_0 \in L^2_\sigma(\Omega) 
\cap L^q(\Omega)^2$, the function $w(t) = A^{-\mu}v(t)$ satisfies
\begin{equation}\label{wbound}
  \|w(t)\|_{L^2}^2 + \int_0^t \|\nabla w(s)\|_{L^2}^2\dd s 
  \,\le\, (1+t)^{\alpha^2 k(\alpha)} \exp\Bigl(K_2(\alpha) \|v_0\|_{L^2}^2 
  + K_3\Bigr) (\|v_0\|_{L^q}^2 + \alpha^2)~,
\end{equation}
for all $t \ge 0$. Moreover $K_2(\alpha)$ and $k(\alpha)$ are 
$\OO(1)$ as $\alpha \to 0$. 
\end{proposition}

\begin{proof} Taking the scalar product of both sides
of \eqref{weq} by $w$, we obtain
\begin{align}\nonumber
  \frac12 \frac{\D}{\D t} \|w(t)\|_{L^2}^2 + \|A^{1/2} w(t)\|_{L^2}^2 
  &+ \alpha \langle F_\mu(u^\chi(t),v(t)),w(t)\rangle +  \alpha
  \langle F_\mu(v(t),u^\chi(t)),w(t)\rangle \\ \label{wfirst} 
  &+ \langle F_\mu(v(t),v(t)),w(t)\rangle 
  \,=\, \alpha \langle A^{-\mu}R^\chi(t),w(t)\rangle~.
\end{align}
It is well known that $\|A^{1/2} w\|_{L^2} = \|\nabla w\|_{L^2}$ for all 
$w \in D(A^{1/2}) = L^2_\sigma(\Omega) \cap H^1_0(\Omega)^2$. To bound the 
other terms, we observe that
\begin{align*}
   |\langle F_\mu(u^\chi,v),w\rangle| \,&=\, |\langle (u^\chi\cdot \nabla) 
   v,A^{-\mu}w\rangle|  \,=\, |\langle (u^\chi\cdot \nabla)A^{-\mu}w 
   ,v\rangle| \\ \,&\le\, \|u^\chi\|_{L^\infty} \|A^{\frac12-\mu}w\|_{L^2} 
  \|v\|_{L^2} \,=\,  \|u^\chi\|_{L^\infty} \|A^{\frac12-\mu}w\|_{L^2} 
  \|A^{\mu}w\|_{L^2} \\ 
  \,&\le\, \|u^\chi\|_{L^\infty} \|A^{1/2}w\|_{L^2} \|w\|_{L^2}~,
\end{align*}
where in the last inequality we used the interpolation inequality 
for fractional powers of $A$. The same argument shows that $|\langle F_\mu
(v,u^\chi),w\rangle| \le \|u^\chi\|_{L^\infty} \|A^{1/2}w\|_{L^2} \|w\|_{L^2}$. 
In a similar way, 
\begin{align*}
  |\langle F_\mu(v,v),w\rangle| \,&=\, |\langle (v\cdot \nabla) 
   v,A^{-\mu}w\rangle|  \,=\, |\langle (v\cdot \nabla)A^{-\mu}w 
   ,v\rangle| \\ \,&\le\, \|v\|_{L^4}^2 \|A^{\frac12-\mu}w\|_{L^2} 
   \,\le\, C \|\nabla v\|_{L^2} \|v\|_{L^2} \|A^{\frac12-\mu}w\|_{L^2} 
   \\ \,&\le\, C \|\nabla v\|_{L^2} \|A^{1/2}w\|_{L^2} \|w\|_{L^2}~.
\end{align*}
Finally, since $|\langle A^{-\mu}R^\chi,w \rangle| =  |\langle R^\chi,
A^{-\mu} w \rangle| \le \kappa_3(1+t)^{-1}\|A^{\frac12-\mu}w\|_{L^2}$ 
by \eqref{Rchi2}, we can use interpolation and Young's inequality
to obtain
\[
  |\alpha\langle A^{-\mu}R^\chi,w \rangle|  \,\le\, \frac{\kappa_3|\alpha|}
  {1+t}\|A^{1/2}w\|_{L^2}^{1-2\mu}\|w\|_{L^2}^{2\mu} \,\le\, \frac14
  \|A^{1/2}w\|_{L^2}^2 + \frac{\|w\|_{L^2}^2}{(1+t)^{\gamma_1}} + 
  \frac{C\alpha^2}{(1+t)^{\gamma_2}}~,
\]
for some $\gamma_1,\gamma_2 > 1$ satisfying $\gamma_2 + 2\mu\gamma_1 = 2$.  
Thus \eqref{wfirst} implies
\begin{align}\nonumber
  \frac{\D}{\D t} \|w\|_{L^2}^2 + \|\nabla w\|_{L^2}^2 \,&\le\, 
  - \|\nabla w\|_{L^2}^2 + C\|\nabla w\|_{L^2} \|w\|_{L^2} (|\alpha| 
  \|u^\chi\|_{L^\infty}  +  \|\nabla v\|_{L^2}) \\ \label{wsecond} 
  &~\quad + \frac12 \|\nabla w\|_{L^2}^2 + \frac{2\|w\|_{L^2}^2}{(1+t)^{\gamma_1}} 
  + \frac{2C\alpha^2}{(1+t)^{\gamma_2}} \\ \nonumber
  \,&\le\, C_1 \|w\|_{L^2}^2 \Bigl(\alpha^2 \|u^\chi\|_{L^\infty}^2 
  + \|\nabla v\|_{L^2}^2 + \frac{1}{(1+t)^{\gamma_1}}\Bigr)
  + \frac{C_2\alpha^2}{(1+t)^{\gamma_2}}~,
\end{align}
for some positive constants $C_1, C_2$.

Now, using \eqref{uchi1} with $p = \infty$ and the logarithmic energy 
estimate \eqref{vlog}, we obtain
\begin{align*}
  C_1\int_0^t \Bigl(\alpha^2 \|u^\chi(s)\|_{L^\infty}^2 
  &+ \|\nabla v(s)\|_{L^2}^2 + \frac{1}{(1+s)^{\gamma_1}}\Bigr)\dd s\\
  \,&\le\, \alpha^2 k(\alpha) \log(1+t) + K_2(\alpha)\|v_0\|_{L^2}^2 
  + C_3~,
\end{align*}
where $K_2(\alpha) = 2^{c|\alpha|}C_1K_1$, $k(\alpha) = C_1a_\infty^2 + 
K_2(\alpha)$, and $C_3 = C_1(\gamma_1-1)^{-1}$. Applying 
Gronwall's lemma to \eqref{wsecond}, we thus find
\[
  \|w(t)\|_{L^2}^2 + \int_0^t \|\nabla w(s)\|_{L^2}^2\dd s 
  \,\le\, (1+t)^{\alpha^2 k(\alpha)} \exp\Bigl(K_2(\alpha) \|v_0\|_{L^2}^2 
  + C_3\Bigr) (\|w_0\|_{L^2}^2 + C_4\alpha^2)~,
\]
where $C_4 = C_2(\gamma_2-1)^{-1}$, and \eqref{wbound} follows since
$\|w_0\|_{L^2} \le C\|v_0\|_{L^q}$ by Lemma~\ref{Amu}. 
\end{proof}

\begin{corollary}\label{vcor}
Under the assumptions of Proposition~\ref{wenergy}, there 
exists a positive constant $K_4$ depending on $|\alpha|$ and 
$\|v_0\|_{L^2 \cap L^q}$ such that, for any $t \ge 2$, there exists a time 
$t_0 \in [t/2,t]$ for which
\begin{equation}\label{vest}
  \|v(t_0)\|_{L^2}^2 \,\le\, K_4 (1+t_0)^{\alpha^2 k(\alpha) - 2\mu}~.
\end{equation}
\end{corollary}

\begin{proof} Fix $t \ge 2$. In view of \eqref{wbound}, 
there exists a time $t_0 \in [t/2,t]$ such that
\[
  \|\nabla w(t_0)\|_{L^2}^2 \,\le\, \frac{2}{t}\int_{t/2}^t 
  \|\nabla w(s)\|_{L^2}^2\dd s \,\le\, \frac{2}{t}
  \,K_0(1+t)^{\alpha^2 k(\alpha)} \,\le\, 2^{\alpha^2k(\alpha)+2}
  K_0(1+t_0)^{\alpha^2 k(\alpha) - 1}~,
\]
where $K_0 = \exp(K_2\|v_0\|_{L^2}^2 + K_3) (\|v_0\|_{L^q}^2 +
\alpha^2)$. Moreover, $\|w(t_0)\|_{L^2}^2 \le K_0(1+t_0)^{\alpha^2
k(\alpha)}$ by \eqref{wbound}. Thus, we obtain \eqref{vest} using 
the interpolation inequality $\|v(t_0)\|_{L^2} = \|A^\mu w(t_0)\|_{L^2} 
\le \|\nabla w(t_0)\|_{L^2}^{2\mu} \,\|w(t_0)\|_{L^2}^{1-2\mu}$.
\end{proof}

We are now able to conclude our sketch of the proof of 
Theorem~\ref{main}. Given $\alpha \in \R$ and $v_0 \in L^2_\sigma(\Omega) 
\cap L^q(\Omega)^2$, let $v(x,t)$ be the solution of the perturbation 
equation \eqref{veq2}. For any $t \ge 2$, we choose $t_0 \in [t/2,t]$ 
as in Corollary~\ref{vcor}, and we apply estimate \eqref{vpoly}. 
We thus obtain
\begin{equation}\label{finest}
  \|v(t)\|_{L^2}^2 \,\le\, 2^{2b_\infty|\alpha|}\Bigl(\|v(t_0)\|_{L^2}^2 + 
  \frac{\kappa_3^2\,\alpha^2}{1+t_0}\Bigr) \,\le\, C_\alpha 
  (1+t)^{\alpha^2 k(\alpha) - 2\mu}~,
\end{equation}
where $C_\alpha > 0$ is $\OO(1)$ as $\alpha \to 0$. Now, if $|\alpha|$
is small enough to that $\alpha^2 k(\alpha) < 2\mu$, the right-hand
side of \eqref{finest} converges to zero (at a suboptimal rate) 
as $t \to \infty$. In particular, the perturbation $v(\cdot,t)$ 
becomes very small in energy norm for large times. In that 
regime, the perturbation equation \eqref{veq2} can be solved 
by a global fixed point argument, which allows to show that
\begin{equation}\label{conv2}
  \|v(\cdot,t)\|_{L^2(\Omega)} + t^{1/2} \|\nabla v(\cdot,t)
  \|_{L^2(\Omega)} \,=\, \OO(t^{-\mu})~,
\end{equation}
as $t \to \infty$, see \cite[Section~3]{GM}. Finally \eqref{conv} 
follows from \eqref{conv2}, because $v(x,t) = u(x,t) - \alpha 
u^\chi(x,t)$ and $\|u^\chi(\cdot,t) - \Theta(\cdot,t)\|_{L^2} + 
\|\nabla u^\chi(\cdot,t) - \nabla \Theta(\cdot,t)\|_{L^2} \le C(1+t)^{-1}$ 
for all $t \ge 0$. \hfill$\Box$

\section{Convergence of the vorticity}\label{sec5}

This final section is devoted to the proof of Proposition~\ref{vort2}. 
We first show that, under the assumptions of Theorem~\ref{main}, one
can control the $L^1$ norm of the vorticity sufficiently far away 
from the obstacle $K$. 

\begin{proposition}\label{noleak}
Under the assumptions of Proposition~\ref{vort2}, the vorticity
$\omega = \curl u$ satisfies
\begin{equation}\label{leakest}
  \int_{|x|\ge t^{1/2}\log t}|\omega(x,t)| \dd x \,=\, \OO(t^{-\mu})~,
  \qquad \hbox{as} \quad t \to \infty~.
\end{equation}
\end{proposition}

\begin{proof}
Since by \eqref{omloc} the initial vorticity is assumed to be 
square integrable, the solution $u(x,t)$ of \eqref{NS} given 
by Theorem~\ref{main} satisfies
\begin{equation}\label{ubound}
  \|u(\cdot,t) -\alpha\Theta(\cdot,t)\|_{L^2(\Omega)} + 
  (1+t)^{1/2} \|\nabla u(\cdot,t) -\alpha\nabla\Theta(\cdot,t)
  \|_{L^2(\Omega)} \,\le\, \frac{C_0}{(1+t)^\mu}~,
\end{equation}
for all $t \ge 0$, where $C_0 > 0$ depends only on the initial
data. The associated vorticity $\omega = \curl u$ is a solution
of the advection-diffusion equation
\begin{equation}\label{omeq}
  \partial_t \omega + u\cdot\nabla\omega \,=\, \Delta \omega~, 
  \qquad x \in \Omega~, \quad t > 0~,
\end{equation}
but the no-slip boundary condition becomes very complicated when
expressed in terms of $\omega$. It is thus difficult to use
\eqref{omeq} to obtain estimates in the whole domain $\Omega$, and in
particular near the boundary $\partial\Omega$. Here, however, our goal
is to bound $\omega$ near infinity, so we can avoid that problem 
using localized energy estimates and invoking \eqref{ubound} to 
control the flux terms in the regions where the localization 
function is not constant.

Given $T \ge 4$ and $R \ge 1$, we define the cut-off function
\begin{equation}\label{psidef}
  \psi(x,t) \,=\, \phi\Big(\frac{|x|}{r(t+T)}\Bigr)\,
  \Bigl(1 - \phi\Bigl(\frac{|x|}{R}\Bigr)\Bigr)~, \qquad 
  x \in \R^2~, \quad t \ge 0~,
\end{equation}
where $r(t) = 2^{-3/2}t^{1/2}\log(t/2)$, and $\phi : [0,\infty) \to
[0,1]$ is a smooth, nondecreasing function satisfying $\phi(r) = 0$
for $r \le 1$ and $\phi(r) = 1$ for $r \ge 2$. We always assume that
$T \ge 4$ is large enough so that the support of $\psi(\cdot,t)$ is
contained in $\Omega$, and that $R \ge 1$ is large enough (depending
on $t$ and $T$) so that $\psi(\cdot,t)$ is not identically zero. Given
$\lambda > 0$, we also denote
\[
  \Phi_\lambda(\omega) \,=\, (\lambda^2 + \omega^2)^{1/2} - \lambda~,
\]
and we observe that $0 \le \Phi_\lambda(\omega) \le |\omega|$ and 
$\Phi_\lambda''(\omega) \ge 0$ for all $\omega \in \R$. 

Now, using \eqref{omeq}, we obtain by a direct calculation
\begin{align*}
  \frac{\D}{\D t}\int_\Omega \psi \Phi_\lambda(\omega) \dd x 
  \,&=\, \int_\Omega \Bigl(\psi_t + \Delta \psi + (u\cdot\nabla)
  \psi\Bigr)\Phi_\lambda(\omega)\dd x - \int_\Omega \psi\Phi_\lambda''(\omega)
  |\nabla \omega|^2 \dd x \\
  \,&\le\, \int_\Omega \Bigl(\Delta \psi + (u\cdot\nabla)
  \psi\Bigr)\Phi_\lambda(\omega)\dd x~,
\end{align*}
because $\psi_t \le 0$ and $\Phi_\lambda''(\omega) \ge 0$. If we integrate 
this inequality over $t \in [0,T]$, we find
\[
  \int_\Omega \psi(x,T) \Phi_\lambda(\omega(x,T)) \dd x \,\le\, 
  \int_\Omega \psi(x,0) \Phi_\lambda(\omega_0(x)) \dd x + 
  \int_0^T \!\!\int_\Omega Q_{\psi,u}(x,t) \Phi_\lambda(\omega(x,t))\dd x 
  \dd t~,
\]
where $Q_{\psi,u} = |\Delta \psi + (u\cdot\nabla)\psi|$. Using
Lebesgue's monotone convergence theorem, we can take the limit
$\lambda \to 0$ in both sides, and we arrive at the simpler estimate
\begin{equation}\label{om1}
  \int_\Omega \psi(x,T) |\omega(x,T)| \dd x \,\le\, 
  \int_\Omega \psi(x,0) |\omega_0(x)| \dd x + 
  \int_0^T \!\!\int_\Omega Q_{\psi,u}(x,t)|\omega(x,t)|\dd x 
  \dd t~.
\end{equation}

Our next task is to take the limit $R \to \infty$ in \eqref{om1}. 
Again, we use the monotone convergence theorem, except in the 
last integral where it does not apply. To treat that term, we 
observe that $Q_{\psi,u}(\cdot,t)$ vanishes identically except 
in the region $D_R \cup D_{r(t+T)}$, where for any $\rho > 0$ 
we denote $D_\rho \,=\, \{x \in \R^2 \,|\,\rho \le |x| \le 2\rho\}$. 
Taking $R > 0$ sufficiently large and using \eqref{ubound}, we easily 
obtain
\[
  \int_{D_R} Q_{\psi,u}|\omega|\dd x \,\le\, \int_{D_R}
  \Bigl(|\Delta \psi| + |u| |\nabla \psi|\Bigr) |\omega|\dd x 
  \,\le\, \frac{C_1}{R}~,
\]
for some $C_1 > 0$ independent of $t \in [0,T]$. The contribution
of the annulus $D_R$ is therefore negligible for large $R$, hence 
taking the limit $R \to \infty$ in \eqref{om1} we arrive at
\[
  \int_\Omega \psi_1(x,T) |\omega(x,T)| \dd x \,\le\, 
  \int_\Omega \psi_1(x,0) |\omega_0(x)| \dd x + 
  \int_0^T \!\!\int_\Omega Q_{\psi_1,u}(x,t)|\omega(x,t)|\dd x 
  \dd t~,
\]
where $\psi_1(x,t) = \phi(|x|/r(t+T))$. In particular $\psi_1(x,T) = 1$ 
for $|x| \ge T^{1/2}\log T$ and $\psi_1(x,0) = 0$ for $|x| \le 
r(T)$, hence the last inequality implies 
\begin{equation}\label{om2}
  \int_{|x| \ge T^{1/2}\log T} |\omega(x,T)| \dd x \,\le\, 
  \int_{|x| \ge r(T)}|\omega_0(x)| \dd x + 
  \int_0^T \!\!\int_\Omega Q_{\psi_1,u}(x,t)|\omega(x,t)|\dd x \dd t\,.
\end{equation}

To conclude the proof of Proposition~\ref{noleak}, it remains to 
estimate both terms in the right-hand side of \eqref{om2}. First, 
using \eqref{omloc} and H\"older's inequality, we easily find
\begin{align}\nonumber
  \int_{|x| \ge r(T)}|\omega_0(x)| \dd x \,&\le\, 
  \left(\int_\Omega (1{+}|x|^2)^m |\omega_0(x)|^2\dd x\right)^{1/2}
  \left(\int_{|x| \ge r(T)} \frac{1}{(1{+}|x|^2)^m} \dd x
  \right)^{1/2} \\ \label{om3}
  \,&\le\, C r(T)^{-(m-1)} \,\le\, C_2 T^{-\mu}~,
\end{align}
for some $C_2 > 0$ independent of $T$. In the last inequality, 
we used the hypothesis $m \ge 2/q = 1 + 2\mu$ and the fact 
that $r(T) \ge CT^{1/2}$. On the other hand, since $\psi_1(x,t)$ is 
given by \eqref{psidef} with $R = \infty$, there exists $C_3 > 0$ 
such that
\[
  |\nabla\psi_1(x,t)| \,\le\, \frac{C_3}{r(t+T)}\one_{D'_t}~, \qquad
  |\Delta\psi_1(x,t)| \,\le\, \frac{C_3}{r(t+T)^2}\one_{D'_t}~, 
\]
where $D'_t =  D_{r(t+T)} = \{x \in \R^2 \,|\, r(t+T) \le |x| \le 
2r(t+T)\}$. It follows that 
\begin{equation}\label{Qbound}
  \int_\Omega Q_{\psi_1,u}(x,t)|\omega(x,t)|\dd x \,\le\, 
  C_3 \int_{D'_t}\left(\frac{1}{r(t+T)^2} + \frac{|u(x,t)|}{r(t+T)}
  \right)|\omega(x,t)|\dd x~.
\end{equation}
But using \eqref{Qbound} and \eqref{ubound} we easily find, 
\begin{align*}
  \int_{D'_t}|\omega(x,t)|\dd x \,&\le\, \int_{D'_t} \Bigl|\omega(x,t) - 
  \alpha \Xi(x,t)\Bigr| \dd x + |\alpha| \int_{D'_t} \Xi(x,t) \dd x \\
  \,&\le\, \meas(D'_t)^{1/2} \|\omega(\cdot,t) - \alpha \Xi(\cdot,t)
  \|_{L^2(\Omega)} + |\alpha| \|\Xi(\cdot,t)\|_{L^1(D'_t)} \\
  \,&\le\, C\,\frac{r(t+T)}{(1+t)^{\mu + 1/2}} + C\,\exp\Bigl(
  -\frac{r(t+T)^2}{4(1+t)}\Bigr) \,\le\, 
  C\,\frac{r(t+T)}{(1+t)^{\mu + 1/2}}~,
\end{align*}
where in the last inequality we used the fact that $r(t+T) 
\ge C(t+T)^{1/2}\log(t+T)$. In a similar way,
\begin{align*}
  \int_{D'_t}|u(x,t)| &|\omega(x,t)|\dd x \,\le\, 
  \int_{D'_t}|u(x,t)| |\omega(x,t) - \alpha\Xi(x,t)|\dd x +  
  |\alpha|\int_{D'_t}|u(x,t)| |\Xi(x,t)|\dd x \\
  \,&\le\, \|u\|_{L^2(D'_t)} \Bigl(\|\omega(\cdot,t) - \alpha\Xi(\cdot,t)
  \|_{L^2(\Omega)} + \|\Xi(\cdot,t)\|_{L^2(D'_t)}\Bigr) 
  \,\le\, \frac{C}{(1+t)^{\mu + 1/2}}~.
\end{align*}
Inserting these estimates in the right-hand side of \eqref{Qbound}, 
we obtain
\begin{equation}\label{om4}
  \int_0^T\!\!\int_\Omega Q_{\psi_1,u}(x,t)|\omega(x,t)|\dd x \dd t
  \,\le\, \int_0^T \frac{C}{r(t+T)\,(1+t)^{\mu + 1/2}} \dd t
  \,\le\, \frac{C_4}{T^\mu}~,
\end{equation}
for some $C_4 > 0$ independent of $T$. Thus, if we combine \eqref{om2}, 
\eqref{om3}, and \eqref{om4}, we conclude that
\[
  \int_{|x| \ge T^{1/2}\log T} |\omega(x,T)| \dd x \,\le\, 
  \frac{C_2 + C_4}{T^\mu}~,
\]
for all sufficiently large $T > 0$, which is the desired result. 
\end{proof}

\medskip
It is now easy to conclude the proof of Proposition~\ref{vort2}. 
For $t > 0$ sufficiently large, we denote $\Omega_t \,=\, \{x \in 
\Omega \,|\, |x|\le t^{1/2}\log t\}$ and we decompose
\[
  \int_\Omega |\omega(x,t) - \alpha \Xi(x,t)|\dd x \,\le\, 
  \int_{\Omega_t} |\omega(x,t) - \alpha \Xi(x,t)|\dd x 
  + \int_{\Omega\setminus\Omega_t} \Bigl(|\omega(x,t)| + |\alpha| \Xi(x,t)
  \Bigr)\dd x~.
\]
The last integral in the right-hand side is controlled using 
Proposition~\ref{noleak} and the explicit expression \eqref{Xidef} 
of $\Xi(x,t)$. To estimate the first integral, we use \eqref{ubound} 
and H\"older's inequality\:
\[
  \int_{\Omega_t} |\omega(x,t) - \alpha \Xi(x,t)|\dd x \,\le\, 
  \sqrt{\pi}\,t^{1/2}\log t \,\|\omega(\cdot,t) - \alpha\Xi(\cdot,t)
  \|_{L^2(\Omega)} \,\le\, C \frac{\log t}{(1+t)^\mu}~.
\]
Summarizing, we find
\[
  \int_\Omega |\omega(x,t) - \alpha \Xi(x,t)|\dd x \,\le\, 
  C \frac{\log t}{(1+t)^\mu}~,
\]
for all sufficiently large $t > 0$. This concludes the proof. \hfill$\Box$

\section{Appendix\: estimates for truncated Oseen 
vortices}\label{sec6}

In this appendix, we collect a few estimates on the truncated 
Oseen vortices \eqref{uchidef} which are used throughout the
paper. We first list a few bounds which follow from \eqref{uchidef} 
and \eqref{omegachi} by rather straightforward calculations, 
see \cite[Lemma~2.1]{GM}. 

\begin{lemma}\label{uchiest} \quad\\
{\bf 1.} For any $p \in (2,\infty]$, there exists a constant
$a_p > 0$ such that
\begin{equation}\label{uchi1}
 \|u^\chi(\cdot,t)\|_{L^p(\R^2)} \,\le\, \frac{a_p}{(1+t)^{\frac12 -
 \frac1p}}~, \qquad t \ge 0~.
\end{equation}
{\bf 2.} For any $p \in (1,\infty]$, there exists a constant
$b_p > 0$ such that
\begin{equation}\label{uchi2}
 \|\nabla u^\chi(\cdot,t)\|_{L^p(\R^2)} \,\le\, \frac{b_p}{(1+t)^{1 
 - \frac1p}}~, \qquad t \ge 0~.
\end{equation}
{\bf 3.} There exists a constant $\kappa_1 > 0$ such that, for all 
$t,s \ge 0$,
\begin{equation}\label{uchi3}
 \|u^\chi(\cdot,t) - u^\chi(\cdot,s)\|_{L^2(\R^2)}^2 \,\le\, 
 \kappa_1\,\Big|\log \frac{1+t}{1+s}\Big|~. 
\end{equation}
{\bf 4.} There exists a constant $\kappa_2 > 0$ such that, for all 
$t,s \ge 0$,
\begin{equation}\label{uchi4}
 \|\nabla u^\chi(\cdot,t) - \nabla u^\chi(\cdot,s)\|_{L^2(\R^2)}^2 
 \,\le\, \kappa_2\,\Big|\frac{1}{1+t} - \frac{1}{1+s}\Big|~. 
\end{equation}
\end{lemma}

\medskip
Since the truncated Oseen vortex is not a solution of the Navier-Stokes
equation, we also need a control on the remainder term $R^\chi
= \Delta u^\chi - \partial_t u^\chi = (\Delta\chi)\Theta + 2
(\nabla\chi\cdot \nabla)\Theta$, which has the explicit expression
\begin{equation}\label{Rchidef}
  R^\chi(x,t) \,=\,  \Theta(x,t) \Delta \chi(x) + 2\frac{x\cdot\nabla
  \chi(x)}{|x|^2}\Bigl(x^\perp \Xi(x,t) - \Theta(x,t)\Bigr)~. 
\end{equation}

\begin{lemma}\label{Rchiest}
There exists a constant $\kappa_3 > 0$ such that, for any 
$p \in [1,\infty]$,
\begin{equation}\label{Rchi1}
  \|R^\chi(\cdot,t)\|_{L^p(\R^2)} \,\le\, \frac{\kappa_3}{1+t}~,  
  \qquad t \ge 0~.
\end{equation}
Moreover, for any vector field $u \in H^1_\loc(\R^2)^2$, we have
\begin{equation}\label{Rchi2}
  \Bigl|\int_{\R^2} R^\chi(x,t)\cdot u(x)\dd x\Bigr| 
   \,\le\, \frac{\kappa_3}{1+t}\,\|\nabla u\|_{L^2(D)}~, 
  \qquad t \ge 0~,
\end{equation}
where $D \subset \Omega$ is a compact annulus containing the support 
of $\nabla\chi$. 
\end{lemma}

\begin{proof}
It is clear from \eqref{Rchidef} that $|R^\chi(x,t)| \le C
(1+t)^{-1}\mathbf{1}_D(x)$ for all $x \in \R^2$ and all $t \ge 0$,  
and \eqref{Rchi1} follows immediately. Moreover, we have 
$R^\chi(x,t) = x^\perp Q^\chi(x,t)$ for some radially symmetric 
scalar function $Q(x,t)$, hence $R^\chi(\cdot,t)$ has zero mean
over the annulus $D$. If $u \in H^1_\loc(\R^2)^2$ and if we denote by
$\bar u$ the average of $u$ over $D$, the Poincar\'e-Wirtinger 
inequality implies
\[
  \Bigl|\int_{\R^2} R^\chi(x,t)\cdot u(x)\dd x\Bigr| \,=\,
  \Bigl|\int_D R^\chi(x,t)\cdot (u(x) - \bar u)\dd x\Bigr| 
  \,\le\, C \|R^\chi(\cdot,t)\|_{L^2(\R^2)} \|\nabla u\|_{L^2(D)}~, 
\]
and using \eqref{Rchi1} with $p = 2$ we obtain \eqref{Rchi2}. 
\end{proof}

\bibliographystyle{plain}


\end{document}